
%
\magnification=\magstep1
\baselineskip=20pt
\font\typc=cmbx10 scaled \magstep1  
\font\type=cmbx10 scaled \magstep3  
\global\def\sectitle#1\par{\bigbreak
  \leftline{\bf #1}
  \nobreak\medskip\vskip-\parskip
  \message{#1}
  \noindent}
\global\def\ssectitle#1\par{\bigbreak\medskip
  \leftline{\typc #1}
  \nobreak\bigskip\vskip-\parskip
  \message{#1}
  \noindent}
\hrule
\bigskip
\centerline{\type
 Homotopy groups of the complements}
\medskip
\centerline{\type to singular
hypersurfaces, II}
\bigskip
\bigskip
\centerline {{\bf by A.Libgober}\footnote *{Supported by NSF grant}}
\centerline {Department of Mathematics}
\centerline {University of Illinois at Chicago}
\centerline {P.O.B.4348, Chicago, Illinois, 60680}
\centerline {e-mail u11377 @UICVM.UIC.EDU}
\bigskip
\hrule
\bigskip
\bigskip
\par {\bf Abstract.} The homotopy group $\pi_{n-k} ({\bf C}^{n+1}-V)$
 where $V$ is a hypersurface with a singular locus of dimension $k$
 and good behavior at infinity is described  using generic pencils.
 This is analogous to the van
 Kampen procedure for finding a fundamental group of
 a plane curve. In addition we use
 a   certain representation generalizing the Burau
 representation of the braid group. A divisibility theorem is proven
 that shows the dependence of
  this homotopy group on the local type of singularities
 and behavior at infinity. Examples are given showing that this group
 depends on certain global data in addition to local data on
 singularities.
\bigskip
 \bigskip

\ssectitle{0.\ Introduction.}

\bigskip
\par The fundamental groups of complements to
algebraic curves in ${\bf CP}^2$
were studied by O.Zariski  almost 60 years
 ago (cf. [Z]). He showed that these groups are affected by the
  type and position
 of  singularities. Zariski
  and Van Kampen (cf.[vK])
 described a general procedure for calculating these
groups in terms of the behavior of the intersection of the curve
with a generic line while one varies this line in a pencil.
For the curves with mild and
 few  singularities, these fundamental groups are abelian. For
example, if $C$ is an irreducible curve, having only singular points
 near which $C$ can be given in some coordinate
  system by the equation $x^2=y^2$,
 then its complement has  an abelian fundamental group
 (cf. [F],[De]). On the other hand,
   one knows that there is an abundance of curves with
 non-abelian fundamental groups of  complements  (for example,
 branching curves of generic projections on a plane of surfaces
 embedded in some projective space). Some
 explicit calcualtions were made by O.Zariski. For example, for the
  curve given by equation $f_2^3+f_3^2=0$, where $f_k(z_0,z_1,z_2)$ is
  a generic form of degre $k$,
  the corresponding fundamental group is $PSL_2(\bf Z)$ (cf [Z]).
\par If one thinks of the high dimensional analog of these results
 one can immediately notice that the class of  fundamental groups
 of the complements to hypersurfaces in a projective space  coincides
  with the class of fundamental groups of  complements to the
curves. Indeed by a Zariski-Lefschetz type theorem   (cf. [H])
the fundamental group
 of the complement to a hypersurface $V$ is the same as the fundamental
 group of the complement to the intersection  $V \cap H$ inside $H$
 for generic plane $H$. In this paper we will show, however,
  that the homotopy group $\pi _{n-k} ({\bf C}^{n+1} - V )$, where $k$
 is the dimension of the singular locus of $V$, exhibits  properties
 rather similar to the properties of the fundamental group
 discovered by O.Zariski.
(By abuse of notation we often will omit specifying the base point
 in homotopy groups except when it can cause
  confusion). Actually, as in the case of curves (cf. [Z],[vK]),
  we study a somewhat more general case of hypersurfaces in affine space
 (which is motivated by a desire to apply these results  to
 the covers
 of ${\bf CP}^{n+1}$ of arbitrary degree branching locus of which
 contains  $V$).
  It appears that the  group $\pi _{n-k} ({\bf CP}^{n+1}- (V \cup H))$
   where $H$ is "a hyperplane at infinity", at least after
tensoring with ${\bf Q}$, has a geometric (rather than a homotopy
theoretic) meaning: it depends on the "local type
and position" of  singularities of  section of $V$ by a generic
linear subspace of ${\bf C}^{n+1}$ of codimension $k$ (cf. example 5.4).
 This group  has  a description in the spirit of geometric topology
  similar to the one given by the van Kampen theorem in which
 the Artin braid group is replaced by a certain generalization.
(cf.sect.2). We allow  $V$ to have a certain type of singular behavior at
infinity. The information on the homotopy groups which we obtain in
 the affine situation  when  $V$ is non-singular at infinity,
as we show below, is equivalent to the information
 on the   homotopy groups
 of complements to projective hypersurfaces. An a priori
construction like the one with generic projections mentioned above in
the case of fundamental groups seems absent in our context (though the
case of weighted homogeneous hypersurfaces discussed in section 1
provides class of hypersurfaces with one singular point
in ${\bf C}^{n+1}$
and non-trivial $\pi_n$). However,  one can, starting from the
 equation of curves with non-trivial fundamental
groups of the complements, construct equations of hypersurfaces having a
non-trivial higher homotopy of their complements. In particular, one
obtains hypersurfaces with the same local data but with a
distinct higher homotopy of the complement.
\par A more detailed  content of this paper is the following. In section
1 we start with a study of the complements to non-singular hypersurfaces
 in ${ \bf C}^{n+1}$ that, in particular, implies that the dimension
$n-k$, where $k$ as above is the dimension of the singular locus of $V$,
 is the lowest in which non-trivial homotopy groups can appear. It also
allows one to reduce the study of
$\pi_{n-k}$ to the case when $V$ has only  isolated singularities.
Moreover, we also derive a relation between
 $\pi _n ({\bf CP}^{n+1}-V)$ and
$\pi_n ({\bf C}^{n+1}-V)$ where $V$ is a hypersurface with isolated
singularities and no singularities at infinity. In section 2 we outline
 a procedure for finding $\pi _n ({\bf C}^{n+1}-V)$ using a generic
 pencil of hyperplane sections of $V$. Such a pencil defines
 the geometric monodromy homomorphism of the fundamental group
  of the space of parameters of the pencil
 corresponding to non singular members of this pencil. The target of
  this homomorphism is the
 fundamental group of the space of certain embeddings of $V \cap H$ into
 $H={\bf C}^{n}$ where $H$ is a generic element of the pencil.
  The latter group has a natural
  linear representaion over ${\bf Z} [t,t^{-1}]$
 . Then the  homotopy group in question  is expressed in terms of the
 geometric monodromy and this linear representation (in some cases
 information on certain degeneration operators will be needed). In the
 case of curves, the  whole procedure coincides with van Kampen's
 (cf. [vK]),
 the group of embeddings is the Artin braid group and the geometric
  monodromy is the braid monodromy  (cf. [Mo]).  The theorem 2.4
 reduces then to van Kampen's theorem. The linear  representation
 mentioned above  is the classical Burau representation and the object
 which the theorem 2.4 calculates is the Alexander module of the curve
 (cf. [L2]).   Section
 3 describes a  neccesary condition for the vanishing of $\pi_n ({\bf C}
 ^{n+1}-V)$. Use of this result shows that in some cases the
 contribution from  the degeneration operators can be omitted. The
 vanishing results here are parallel to
  the results on  the commutativity  of the
 fundamental groups of  complements in the case of curves. The issue
 of explicit numerical conditions on singularities, which will assure
 the vanishing of the homotopy groups, is more algebro-geometric
  in nature than most issues
 treated here and will be discussed elsewhere. Note also that another
 vanishing result for $\pi_2$ of the complement to an image of a
 generic projection was obtained in [L3]. The next section gives
 restrictions on $\pi_n ({\bf C}^{n+1}-V)$ imposed by the local type
 of  singularities and the behavior of $V$ at infinity. These results
 are generalizations of divisibility theorems for Alexander
 polynomials in
    [L2]. As a corollary, we show that in the case of the absence of
 singularities at infinity the order of $\pi_n ({\bf C}^{n+1}-V)\otimes
{\bf Q} $ coincides with the characteristic polynomial of the monodromy
 operator acting on $H_n$ of the Milnor fibre
 of  (non-isolated)  singularity of the defining equation of $V$
 at the origin in  ${\bf C}^{n+2}$. These
 characteristic polynomials were also considered in [Di]. In the last
 section we give two methods for constructing  hypersurfaces for which
 $\pi_n ({\bf C}^{n+1}-V) \ne 0$ and we calculate the homotopy groups
 in these cases.  The
 first method is based on  a generalization of Zariski's example
 mentioned in the  first paragraph . We show that if $f_k$ denotes
 a generic form of degree $k$ of $n+2$ variables, $p_i$ $i=1,..,n+1$
 are positive integers and $q_i=(\prod_{i=1}^{n+1}  p_i)/p_i$,
  then the order
 $\pi_n ({\bf CP}^{n+1}-(V \cup H)) \otimes {\bf Q}$ ($H$ is a generic
 hyperplane) as a module over ${\bf Q} [t,t^{-1}]$ (cf. sect.1)
 where $V$ is given by equation
 $$f^{p_1}_{q_1}+...f^{p_{n+1}}_{q_{n+1}}=0 \eqno (0.1)$$
  is the characteristic polynomial of the monodromy of the singularity
 $x_1^{p_1}+...x_{n+1}^{p_{n+1}}=0$. The hypersurface (0.1) and the
 non-vanishing of $\pi_n$ was described in [L1] as a consequence of
 the fact that singularities of $(0.1)$ form a finite set in ${\bf CP}
 ^{n+1}$ which is a {\it complete intersection} of ${n+1}$
 hypersurfaces $f_{q_1}=...f_{q_{n+1}}=0$ in ${\bf CP}^{n+1}$.
 The second method, based on the use of Thom-Sebastiani theorem results
 in examples of hypersurfaces with the same collection of singularities
  but with distinct $\pi_n$'s  of the complement.
\par Part of the results in this paper were outlined in the announcement
 [L1]. On the other,  the latter contained number of
  results on the relationship
 of the homotopy groups in question with the Hodge theory of the cyclic
 covers of ${\bf CP}^{n+1}$ branched over $V$; this will be addressed
 elsewhere. Note that since [L1] appeared two more publications
  ([Di] and [Deg] related to the case
   of hypersurfaces in ${\bf CP}^{n+1} $
 came out. Finally, I want to thank Prof.P.Deligne for discussions and
 the Institute for Advanced Study, where a part of this paper was
  written, for support and hospitality.
  \bigskip
\bigskip

\ssectitle{1\ Preliminaries}

\bigskip
\par This section describes the topology of complements to
 non-singular hypersurfaces ${\bf C}^{n+1}$, the case of weighted
 homogeneous hypersurfaces, the homology of the complements
  and the relationship between complements in affine and projective
 spaces. Note that the complements to nonsingular hypersurfaces
 in projective case were first studied in [KW].
\bigskip
\par {\bf Lemma 1.1.} Let $V$ be a non-singular hypersurface
 in ${\bf C}^{n+1}$
 which is transversal to the hyperplane at infinity (resp. a non-
singular hypersurface in ${\bf CP}^{n+1}$). Then $\pi _i ({\bf
 C}^{n+1}-V)=0$ for $1 <i \le n$ and $\pi_1({\bf C}^{n+1}-V)={\bf Z}$
 (resp. ${\bf Z} / d {\bf Z}$, $d$ is the degree of $V$ ).
\bigskip
\par {\bf Proof.} The statement about the fundamental groups follows
immediately
from the Zariski theorem by taking a section of the hypersurface by
generic plane. A hypersurface satisfying the conditions
of the lemma is isotopic to the hypersurface given by equation:
$z_1^d+...z_{n+1}^{d}=1.$  (resp. projective closure of this). The
$d$-fold cover of ${\bf CP}^{n+1}$
branched over the projective closure of the latter can be identified
with the hypersurface $\cal V$ in ${\bf C}^{n+2}$ given by
$z_0^d+z_1^d+...z_{n+1}^d=1$. This is $n$-connected because, for example,
 it is diffeomorphic
to the Milnor fibre of the isolated singularity $z_0^d+...z_{n+1}^d=0$.
Hence our claim follows in the projective case. The $d$-fold cover of
the
 complement in ${\bf C}^{n+1}$ is obtained by removing from $\cal V$ the
hyperplane section $H_0$ given by $z_0=0$. In other words this $d$-fold
cover is the
hypersurface in ${\bf CP}^{n+1}-(H_0 \cap H_\infty) ={\bf C}^*$
( $H_\infty$ is the hyperplane at infinity) which is  transversal to
these
hyperplanes. Hence the Lefschetz type theorem implies that
$\pi_i ({\cal V}-H_0)=\pi _i({\bf C}^*)$.
\bigskip
\par {\bf Corollary.1.2.} If $V$ is a non-singular hypersurface
transversal
 to the hyperplane at infinity then ${\bf C}^{n+1} - V$ is homotopy
 equivalent to the wedge of spheres $S^1 \vee S^{n+1} \vee... \vee
 S^{n+1}$.
\bigskip
\par {\bf Proof}. Indeed the lemma implies that CW-complex
 ${\bf C}^{n+1}$ is a $({\bf Z},n+1)$ complex in the sense of [Dy]
 and hence is homotopy  equivalent to the wedge as above (i.e. as
 a consequence of stable triviality of $\pi_{n+1} ({\bf C}^{n+1}-V)$
 (cf. [Wh],th 14) and the fact that
    stably trivial modules over ${\bf Z} [t,t^{-1}]$ are free.).
\bigskip
\par {\bf Remark 1.3.} It is interesting to see how such wedge comes up
 geometrically. Let us compare the complement to
$q(z_1,...,z_{n+1})=z_1^2 +...z_{n+1}^2=1$ in ${\bf C}^{n+1}$ with the
complement to $z_0^2+...z_{n+1}^2=0$. The last complement fibres over
 ${\bf C}^*$
 using the map $(z_1,...,z_{n+1}) \rightarrow z_0^2+...+z_{n+1}^2$.
 The fibre is homotopy equivalent to $S^n$. Hence the complement to
the singular quadric can be identified with $S^1 \times S^n$. On the
other hand the degeneration of the non-singular quadric into the
singular one results
in the collapse of the vanishing cycle $S^n$ which is the boundary of
a relative vanishing cycle. This relative vanishing cycle
 can be given explicitly as the set
 $\Delta$ of $(z_1,...,z_{n+1}) \in {\bf R}^{n+1} \subset {\bf C}^{n+1}
 , \vert z_1 \vert ^2+...\vert z_{n+1}
\vert ^2 \le 1-\epsilon$). The complement to $z_0^2+...z_{n+1}^2=1$
 hence can be obtained from the complement to
 $z_0^2+...z_{n+1}^2=0$ by attaching a $(n+1)$-cell. Therefore  the
complement to the non-singular quadric can be identified with  $S^1
\times S^n \cup _{* \times S^n} e_{n+1}= S^1 \vee S^{n+1}$.
\bigskip
 \par {\bf Definition 1.4.} Let $V$ be a hypersurface in ${\bf
CP}^{n+1}$
 and $H$ be a hyperplane. A point of $V$ will be called a singular
 point at infinity if it is a singular point of $V \cap H$. The
 subvariety of singular points of $V$ will be denoted
$Sing_{\infty}(V)$.
\bigskip
\par {\bf Lemma 1.5.}  Let $V$ be a hypersurface in ${\bf CP}^{n+1}$
having the dimension of
 $Sing(V) \cup Sing_{\infty}(V)$  (resp. $ Sing (V)$)
 equal to $k$. If ${\bf C}^{n-k+1}$ is a generic linear subspace of
 ${\bf C}^{n+1}$  of codimension
 $k$ then ${\bf C}^{n-k+1} \cap V$ has isolated singularities and
 $\pi_{n-k} ({\bf C}^{n+1}-V)=\pi_{n-k} ({\bf C}^{n-k+1}-V \cap
 {\bf C}^{n-k+1})$ (resp. $\pi _{n-k} ({\bf CP}^{n+1}- (V \cup H))=
 \pi_{n-k} ({\bf CP}^{n-k+1}-(V \cup H) \cap {\bf CP}^{n-k+1})$).
 Moreover $\pi _1 ({\bf CP}^{n+1}-V \cup H)= {\bf Z}$
 (resp. $\pi _1({\bf CP}^{n+1}-V)={\bf Z} / d{\bf Z}$) and
 $\pi_i ({\bf CP}^{n+1}-V \cup H)=\pi _i({\bf CP}^{n+1}-V)=0)$ for
 $2 \le i <n-k$.
\bigskip
\par {\bf Proof.} The first part is a consequence of the Lefschetz
 theorem (cf. [H]). If $L$ is a generic subspace of codimension $k+1$ in
${\bf CP}^{n+1}$ then $V \cap L$ is a non-singular hypersurface in $L$
 which is transversal to $L \cap H$ and the claim follows from the lemma
 1.1.
  \bigskip
\par {\bf Lemma 1.6.} Let $V$ be a hypersurface with isolated
 singularities including singularities at infinity.
  Then $H_i ({\bf CP}^{n+1}-(V \cup
H),{\bf Z})$ =0 for $2 \le i \le {n-1}$. If $V$ is non-singular then the
 vanishing also takes place for $i=n$.
 If $n \ge 2$ then for $i=n$ this group is isomorphic to
 $H^{n+1}(V,H \cap V,{\bf Z})$.
Moreover $H_1 ({\bf CP}^{n+1}-(V \cup H),{\bf Z})$ is isomorphic to
 ${\bf Z}$ unless $dim V=1$,
in which case this group is the free abelain group of the rank
 equal to the number of irreducible components of $V$.
\bigskip
\par {\bf Proof}. The calculation of the group in the lemma for
 $i \le n-1$ using the Lefschetz theorem can be reduced to the case
$i=n$  (as was done above for homotopy groups). For $i=n$
, as follows from the exact sequence of the pair $({\bf CP}^{n+1}-H,
 {\bf CP}^{n+1}- (V \cup H))$ the group in the
lemma is isomorphic to $H_{n+1} ({\bf CP}^{n+1}-H,{\bf CP}^{n+1}-
 (V \cup H),{\bf Z})$. If $T(V)$ and $T(H)$ are
  the regular neighbourhoods of $V$ and $H$ respectively
  in ${\bf CP}^{n+1}$ and $\partial _\infty$ is intersection
 of the boundary of  $T(H)$  with $T(V)-T(V) \cap T(H)$
then the last homology group can be replaced using excision and
 the Lefschetz duality by $H^{n+1} (T(V), \partial _\infty,{\bf Z})$.
 If $V$ is non-singular then the latter group can be identified
 with $H_{n-1} (V-V \cap H,{\bf Z})=0$, which proves the lemma for
$2 \le i \le n-1$. In the general case the last cohomology group by
 excision is isomorphic to $H^{n+1} (V,H \cap V,{\bf Z})$ which proves
 the lemma. The remaining case of curves is well known (cf. [L]).
\bigskip
\par {\bf Lemma 1.7.} If $V$ and $V \cap H$ are ${\bf Q}$-manifold then
 $H^{n+1} (V,H \cap V,{\bf Q})=0$.
\bigskip
\par {\bf Proof.} This is equivalent to injectivity of
$H^{i} (V,{\bf Q}) \rightarrow  H^{i} (H, {\bf Q})$  for $i=n+2$
 and surjectivity for $i=n+1$.    This follows from
 the Poincare duality with ${\bf Q}$-coefficients and the Lefschetz
theorem
 (the dual homology groups are isomorphic to ${\bf Q}$  or to $0$
 depending on parity of $n$ ).
\bigskip
\par {\bf Remark 1.8.} An easily verifyable condition for a hypersurface
 $V$ to be a ${\bf Q}$-manifold is the following: if $V$ has only
 isolated singularities and for each of the singularities
  the characteristic
 polynomial of the monodromy operator does not vanish at $1$ then $V$
is a $\bf Q$-manifold. This is an immediate consequence of the fact
that the condition on the monodromy is equivalent to the condition that
the link of each singualarity is a ${\bf Q}$ -sphere and of the Zeeman
 spectral sequence (cf. [McC]).
\bigskip
\par Let $V$ be a hypersurface in ${\bf C}^{n+1}$ having only isolated
 singularities including infinity. According to lemma 1.5, one can
 identify $\pi_n ({\bf C}^{n+1}-V)$ with $H_n (\widetilde {
 {\bf C}^{n+1}-V},{\bf Z})$ where $\widetilde {{\bf C}^{n+1}-V}$ is
 the universal cover of the space in question. The group ${\bf Z}$ of
 deck transformations acts on $H_n (\widetilde {{\bf C}^{n+1}-V},
 {\bf Z})$. This action on $\pi_n ({\bf C}^{n+1}-V)$ can be described
 as $\beta
 \rightarrow [\alpha , \beta]- \beta$ where $\beta \in \pi_n ({\bf C}^
{n+1}-V)$, $\alpha \in \pi_1 ({\bf C}^{n+1}-V)$ and $[,]$ is the
Whitehead product (alternatively this action is the one given by
the change of the base point ). The structure
 of $\pi _n ({\bf C}^{n+1}-V)$ as the module over the group ring of the
 fundamental group, i.e. over ${\bf Z} [t,t^{-1}]$, becomes particularly
   simple after tensoring with ${\bf Q}$:
$$\pi _n ({\bf C}^{n+1}-V) \otimes {\bf Q} = \oplus _i {\bf Q} [t,t^{-1}]
 /(\lambda_i) \eqno (1.1)$$
 as modules over ${\bf Q} [t,t^{-1}]$ for some polynomials $\lambda_i$
defined up to a unit of ${\bf Q} [t,t^{-1}]$.
\bigskip
\par {\bf Definition 1.9.} The product $\prod _i \lambda _i$
 is called the order of $\pi _n ({\bf C}^{n+1}-V) \otimes {\bf Q}$
 (as a module over ${\bf Q} [t,t^{-1}]$).
\bigskip
\par {\bf Remark 1.10.} In the low dimensional cases, if one works with
 the homology of infinite cyclic covers one obtains results similar
 to what follows. The order of the corresponding ${\bf Q} [t,t^{-1}]$
 module in the case $n=1$ is the Alexander  polynomial of the curves
 studied in [L2]. Note that if $n=0$ then the one-dimensional homology
 over ${\bf Z}$
  of the infinite cyclic cover of ${\bf C}-V$ (i.e. of the complement
 to, say,  $d$ points) is a free module over ${\bf Z} [t,t^{-1}]$ of
 rank $d-1$.
\bigskip
\par {\bf Lemma 1.11.} Let $f(z_1,...,z_n)$
be a weighted homogeneous polynomial
having an isolated singularity at the origin and $V$ is given by $f=0$ .
Then  $\pi_n ({\bf C}^{n+1}- V)=H_n (M_f,{\bf Z})$ where $M_f$ is the
Milnor fibre of the singularity of $f$. This isomorphism is
an isomorphism of ${\bf Z}[t,t^{-1}]$-modules where the structure
 of such module on $H_n(M_f,{\bf Z})$ is given $t$ acting
 as the monodromy operator.
\bigskip
\par {\bf Proof.} ${\bf C}^{n+1}-V$ can be retracted on the complement
of the link of the singularity of $f$. Hence this follows from the exact
 sequence of fibration and $n-1$-connectedness of the Milnor fibre.
\bigskip
\par {\bf Lemma 1.12.} Let $P_V(t)$ be the order of
 $\pi _n ({\bf CP}^{n+1}
 -(V \cup H)) \otimes {\bf Q}$ as ${\bf Q} [t,t^{-1}]$-module. If
 $H^{n+1} (V,H \cap V,{\bf Q})=0$ then $P_V (1) \ne 0$. In particular
 this homotopy group is a torsion module.
\bigskip
\par {\bf Proof.} Let  us consider the exact sequence corresponding
  to the following exact sequence of the chain complexes of the
 universal cyclic cover $\widetilde {({\bf CP}^{n+1}- (V \cup
H))}_\infty$:
 $$ 0 \rightarrow C_* (\widetilde {({\bf CP}^{n+1}-(V \cup H))}_\infty
 ) \rightarrow
  C_* (\widetilde {({\bf CP}^{n+1}- (V \cup H))}_\infty) \rightarrow
 C_* ( {\bf CP}^{n+1}-(V \cup H)) \rightarrow 0 \eqno (1.2)$$
 where the left homomorphism is the map of free ${\bf Q} [t,t^{-1}]$-
 modules induced by the multiplication by $t-1$.
 We obtain:
  $$ \rightarrow H_n (\widetilde {({\bf CP}^{n+1}- (V \cup H))} _\infty
  , {\bf Q})
  \rightarrow H_n (\widetilde {({\bf CP}^{n+1}- (V \cup H))} _\infty
  , {\bf Q}) \rightarrow H_n ({\bf CP}^{n+1}-(V \cup H),{\bf Q})
   \rightarrow
 \eqno (1.3) $$
 (the left homomorphism is the multiplication by $t-1$.) The right group
 in $(1.3)$ is trivial by the assumption and lemma 1.6. Hence
 the multiplication by $t-1$ in
 $\pi _n ({\bf CP}^{n+1}- (V \cup H)) \otimes {\bf Q}=
 H_n (\widetilde {{\bf CP}^{n+1}- (V \cup H)}\infty
 ,{\bf Q})$ is surjective.  Therefore,
 its cyclic decomposition does not have either free summands or
 summands of the form
 ${\bf Q} [t,t^{-1}]/(t-1)^{\kappa} {\bf Q}[t,t^{-1}]$ $\kappa
\in {\bf N}$.
\bigskip
\par {\bf Lemma 1.13.} Let $ H$ be a generic hyperplane and $V$ a
 hypersurface of degree $d$ with isolated singularities in
 ${\bf CP}^{n+1}$.
  Let $d {\bf Z}$ be the subgroup of ${\bf Z}$ of index $d$.
  Then $\pi_n ({\bf CP}^{n+1}- V)$ is isomorphic to the covariantes
 of $\pi _n ({\bf CP}^{n+1}- (V \cup H))^{d {\bf Z}}$ (i.e. the quotient
 by the submodule of images by the action of augmentation ideal of the
 subgroup: $\pi_n /
 (t^d-1) \pi_n$) with the standard  action of ${\bf Z} / d {\bf Z}$.
\bigskip
\par {\bf Proof}. First let us show  that the module of covariants
in the lemma is isomorphic to  $$H_n (\widetilde {({\bf CP}^{n+1}-
(V \cup H) )_d} ,{\bf Z}) \eqno (1.4)$$   where
$\widetilde {({\bf CP}^{n+1}-  (V \cup H))_d}$
 is the $d$-fold cyclic cover of the corresponding
 space. Indeed the sequence of the chain complexes similar to (1,2):
 $$0 \rightarrow C_* ( \widetilde {({\bf CP}^{n+1}- (V \cup H))}_
\infty)
 \rightarrow C_* ( \widetilde {({\bf CP}^{n+1} -(V \cup H))}_ \infty)
 \rightarrow  C_*  (\widetilde {({\bf CP}^{n+1}- (V \cup H))}_d)
 \rightarrow  \eqno (1.5)$$
 where the left homomorphism is the multiplication by $t^d-1$,
 gives rise to  the homology sequence:
 $$\rightarrow H_n (\widetilde {({\bf CP}^{n+1}- (V \cup H))}_\infty
,{\bf Z})
 \rightarrow H_n (\widetilde {({\bf CP}^{n+1}- (V \cup H))} _\infty,
{\bf Z})
 \rightarrow H_n (\widetilde {({\bf CP}^{n+1}- (V \cup H))} _d, {\bf Z})
 \rightarrow \eqno (1.6)$$  in which the
 right homomorphism is surjective because of the vanishing of
 $\pi_ {n-1} ({\bf CP}^{n+1}- (V \cup H))$, which proves our claim.
 \par To conclude the proof of the lemma we need to show that
 $$H_n (\widetilde {({\bf CP}^{n+1}- (V \cup H))}_d,{\bf Z})=
   H_n (\widetilde {({\bf CP}^{n+1}-V)}_d,{\bf Z}) \eqno (1.7)$$ where
$\widetilde  {({\bf CP}^{n+1}-V)}_d$ is the universal
 cyclic cover of ${\bf CP}^{n+1}-V$ . This will
follow from  the vanishing of the relative group $H_i ({\widetilde ({\bf
CP}^{n+1}-V)}_d,$ $({\widetilde ({\bf CP}^{n+1}-(V \cup H))}_d ,{\bf
Z})$  for $i=n,n+1$. Let $W_d$ be the cyclic $d$-fold branched over $V$
covering of ${\bf CP}^{n+1}$ and let $Z \subset W_d$ be the
ramification locus (isomorphic to $V$). Let $H_d$ be the submanifold of
$W_d$ which maps onto $H \cap V$. Let us consider the regular
neighbourhoods $T(H_d)$ and $T(Z)$ in $W_d$ of $H_d$ and $Z$
respectively. The boundary of $T(H_d)-T(H_d) \cap T(Z)$ contains the
part $\partial _1$ which is the part of the boundary of $T(H_d)$.
The complementary to $\partial _1$ part of the boundary of
$T(H_d)-T(H_d) \cap T(Z)$ we denote $\partial _2$.  By excision:
$$H_i (\widetilde {({\bf CP}^{n+1}-V)}_d, \widetilde {({\bf CP}^{n+1}-
(V \cup H)}_d,{\bf Z})=H_i (T(H_d)-T(H_d) \cap T(Z),\partial _1,{\bf
Z}) \eqno (1.8)$$
 Moreover $H_i (T(H_d)-T(H_d) \cap T(Z),\partial _1,{\bf Z})=$
$H^{2n+2-i} (T(H_d)-T(H_d) \cap T(Z), \partial _2,{\bf Z})=$
$H^{2n+2-i} (H_d,H_d \cap Z,{\bf Z})$ by duality, excision and
retraction. The assumption that $H$ is generic implies that $H_d$,which
is a cyclic branched cover of $H$ of degree $d$ with branching locus $H
\cap Z$, is non-singular and hence
$H^{2n+2-i} (H_d,H_d \cap Z,{\bf Z})= H_{i-2} (H_d - H_d \cap Z,{\bf
Z})$. The latter group which is the homology group of affine
hypersurface
transversal to the hyperplane at infinity is trivial except for
$i=2$ and $i=n+2$. This implies the lemma.
\bigskip
\bigskip
\ssectitle{2.\ Calculation of the homotopy groups using generic
pencils.}

\bigskip
\par In this section we shall describe a method for calculation of
 $\pi_n ({\bf C}^{n+1}-V)$ using the monodromy action on the homotopy
groups
 of the complement in a generic element of a linear pencil of
hyperplanes to the intersection of this element with $V$ obtained
 by moving this element around a loop in the parameter space of the
 pencil.
 Monodromy action is the  composition of "geometric monodromy" with
 values in the fundamental group of certain embeddings and a linear
 representation of this group over ${\bf Z} [t,t^{-1}]$.
In the case  $n=1$ this construction with monodromy taking values in
 the fundamental group of the space of embeddings (which in this
 case is the Artin's braid group), reduces to the van Kampen theorem
(cf.
 [vK]). The composition of this monodromy with the Burau representation
 leads to calculation of the Alexander module of the curve.
 We will start with specifying the loops such that by moving
along these
 loops we  will get the information needed.
 \bigskip
\par  {\bf Definition 2.1.} Let $t_1,...,t_N$ be a finite set of points
in ${\bf C}$. A system of generators of $\gamma_i \in \pi_1 ({\bf C}-
\bigcup _i t_i, t_0)$ is called {\it good} if each of the loops
 $\gamma_i: S^1
\rightarrow  {\bf C}- \bigcup_i t_i$ extends to a map of the disk
$D^2 \rightarrow  {\bf C}$ with non-intersecting images for distinct
$i$'s.
\bigskip
 \par A standard method for constructing a good system of generators is
to select a system of small disks $\Delta _i$ about each of $t_i$
$i=1,...,N$, to
choose a system of $N$ non-intersecting paths $\delta_i$ connecting the
base point $t_0$ with a point of $\partial \Delta _i$  and to take
$\gamma_i = \delta ^{-1} \circ \partial \Delta _i \circ \delta_i$
(with, say, the counterclockwise orientation of $\partial \Delta _i$).
 \par Let $V$ be a non-singular hypersurface in ${\bf C}^{n}$ which
 is transversal to the hyperplane at infinity. Let us consider a sphere
 in $S^{2n-1}$ in ${\bf C}^{n}$ of a sufficiently large radius. Let
  $\partial_ \infty V= V \cap S^{2n-1}$. Let us consider the space of
   $Emb (V,{\bf C}^{n})$ of submanifolds of ${\bf C}^n$ which are
   diffeomorphic to $V$ and which are isotopic to the chosen embedding
 of $V$,  such that
 for any $V^{\prime} \in Emb (V,{\bf C}^{n})$ one has $V^{\prime}
  (V) \cup S^{2n-1}=\partial _\infty V$.
   We assume the compact open topology
 on this space of submanifolds. Let us  describe a certain linear
 representation
 of $\pi_1 (Emb (V,{\bf C}^{n}))$
  (over ${\bf Z} [t,t^{-1}]$) which after a choice
 of a basic gives homomorphism into $GL_r ({\bf Z} [t,t^{-1}])$ where
 $r$ is the rank of $\tilde H_n ({\bf C}^{n}-V,{\bf Z})$ (the reduced
 homology of the complement).
\par  Let $Diff ({\bf C}^{n},S^{2n-1})$ be the
group of diffeomorphisms of ${\bf C}^{n}$ which act as the identity
 outside of $S^{2n-1}$. This group can be identified with $Diff
  (S^{2n},D_{2n})$ of the diffeomorphisms of the sphere fixing
 a disk (cf. [ABK]). Let $Diff ({\bf C}^{n},V)$ be the subgroup
 of $Diff ({\bf C}^{n},S^{2n})$  of the diffeomorphisms  which
 take $V$ into itself. The group $Diff ({\bf C}^{n},S^{2n-1})$ acts
 transitively on $Emb (V,{\bf CP}^{n})$ with the stabilizer $Diff ({\bf
 C}^{n},V)$ which implies the following exact sequence:
 $$\pi_1 (Diff ( S^{2n+2},D_{2n+2})) \rightarrow  \pi_1 (Emb ({\bf C}
 ^{n+1},V)) \rightarrow \pi_0 (Diff ({\bf C}^{n+1},V)) \rightarrow
\pi_0  (Diff( S^{2n+2},D^{2n+2})) \rightarrow \eqno (2.1)$$
 Any element in $Diff ({\bf C}^{n},V)$ induces the self map  of
 ${\bf C}^{n}-V$ and the self map of the universal (cyclic in the case
  $n=1$ )
  cover of this space. Hence it induces an automorphism of $H_n
(\widetilde {{\bf  C}^n-V,{\bf Z}})=\pi_n ({\bf C}^n-V)$, $n>1$.
.  The
composition of the boundary
 homomorphism in (2.1) with the map of $\pi _0 (Diff ({\bf C}^n,V))$
just described results in: $$\lambda: \pi_1 (Emb ({\bf C}^n,V))
 \rightarrow  Aut (\pi_n ({\bf C}^n-V)) \eqno (2.2)$$
\par  In the case $n=1$, $V$ is just a collection of points in ${\bf
C}$,
 $\pi_1 (Emb ({\bf C},V))=\pi_0 (Diff ({\bf C},V))$ is  Artin's
 braid group, and this construction gives the homomorphism of the braid
 group into $Aut ( H_1 (\widetilde {{\bf C}-V,{\bf Z}})$ which, after a
choice  of the basis in $ H_1 (\widetilde {{\bf C}-V})$
 corresponding to the choice
 of the generators of the braid group, gives the reduced Burau
representation. This construction coincides
 with the one described in [A].
\par Now we  can define the relevant monodromy operator corresponding to
a  loop in the parameter space of a linear pencil of hyperplane sections.
  Let $V$ be a hypersurface in ${\bf CP}^{n+1}$
 which has only isolated singularities and $H$ be the hyperplane at
 infinity (which we shall assume transversal to $V$).
  Let $L_t$, $t \in {\bf C}$, be
 a pencil of hyperplanes the projective closure of which has the base
locus $B \subset H$ such that $B$ is transversal to $V$. Let $t_1,..
.,t_N$ denote those $t$ for which $V \cap L_t$ has a singularity.
We assume that for any $i$ the singularity of $V \cap L_{t_i}$ is
outside of $H$. Over ${\bf C}- \bigcup_i t_i$ the pencil   $L_t$
 defines a locally trivial fibration $\tau$ of ${\bf C}^{n+1}-V$
 with a non-singular hypersurface in ${\bf C}^n$ as a fibre transversal
to the hyperplane at infinity. The restriction of this  fibration
 on the complement to a sufficiently large ball is trivial, as follows
 from the assumptions on the singularities at infinity.
Let $\gamma: [0,1] \rightarrow  {\bf C}- \bigcup_i t_i$ $(i=1,...,N)$
be a loop with the base point $t_0$.
   A choice of a trivialization of the pull back  of the
fibration $\tau$ on $[0,1]$ using $\gamma$, defines a loop $e_{\gamma}$
in $Emb (L_{t_0},V \cap L_{t_0})$. Different trivializations
 define homotopic loops in this space.
 \par {\bf Definition 2.2.} The monodromy operator corresponding
 $\gamma$ is the  element in $Aut (\pi_n (L_{t_0}-L_{t_0} \cap V))$
 corresponding in (2.2) to  $e_\gamma$.
\par Next we will need to associate the homomorphism with a singular
fibre $L_{t_i}$
 and a loop $\gamma$ with the  base point $t_0$
 in the parameter space of the pencil which there bounds
 a disk $\Delta_{t_i}$ not containing other singular points
  of the pencil :
   $$\pi_{n-1} (L_{t_i}-L_{t_i} \cap V) \rightarrow
 \pi _n (L_{t_0}-L_{t_0} \cap V)/ Im (\Gamma -I) \eqno (2.3)$$
 where $\Gamma$ is the monodromy operator corresponding to $\gamma$.
 \par First let us note that the module on the right in (2.3) is
 isomorphic to the homology $H_n (\widetilde
  {\tau ^{-1} (\partial \Delta _{t_i})} ,{\bf Z})$ of the infinite
 cyclic cover of the restriction of the fibration $\tau$ on the
 boundary of $\Delta _{t_i}$. This follows  immediately  from the Wang
sequence
 of a fibration over a circle and the vanishing of the homotopy of
 $L_{t_0}-L_{t_0} \cap V$ in dimensions below $n$. Let $B_i$ be a
 polydisk in ${\bf C}^{n+1}$ such that $B_i=\Delta _i \times B$ for
 a certain polydisk $B$ in $L_{t_0}$. Then $\widetilde {\tau ^{-1}
  (\Delta _i)-B_i}$ is a trivial fibration over $\Delta_i$ with the
 infinite cyclic cover $\widetilde {L_{t_i}-L_{t_i} \cap V}$ as a fibre.
 In particular,  one obtains the map:
  $$\pi_ {n-1} (L_{t_0}-L_{t_0} \cap V)=
 H_{n-1} (\widetilde {L_{t_0}-L_{t_0} \cap V},{\bf Z}) \rightarrow
 H_n (\widetilde {\tau ^{-1} (\Delta _i)-B_i},{\bf Z})=$$
 $$H_{n-1} (\widetilde{L_{t_0}-L_{t_0} \cap V},{\bf Z}) \oplus
 H_n (\widetilde {L_{t_0}-L_{t_0} \cap V},{\bf Z}) \eqno (2.4)$$
 \par {\bf Definition 2.3.} The degeneration operator is the map (2.3)
 given by composition of the map (2.4) with  the map
 $H_n (\widetilde {\tau ^{-1}
   (\Delta _i)-B_i},{\bf Z}) \rightarrow H_n(\widetilde {\tau^{-1}
 (\Delta _i)},{\bf Z})=\pi _n (L_{t_0}-L_{t_0} \cap V)$ induced by
 inclusion.
\bigskip
\par {\bf Theorem 2.4.} Let $V$ be a hypersurface in ${\bf CP}^{n+1}$
having only isolated singularities and transversal to the hyperplane $H$
at infinity. Consider a pencil of hyperplanes in ${\bf CP}^{n+1}$
the base locus of which belongs to $H$ and is transversal in $H$ to
$V \cap H$. Let
${\bf C}^n_t$ ($t \in \Omega$) be the  pencil of hyperplanes in ${\bf
C}^{n+1}={\bf
CP}^{n+1}-H$ defined by  $L_t$ (where $\Omega={\bf C}$ is the set
parametrizing all elements of the pencil $L_t$ excluding $H$). Denote by
$t_1,...,t_N$  the collection of
those $t$ for which $V \cap L_t$ has a singularity.
We shall assume that the pencil was chosen so that $L_t \cap H$ has at
 most one singular
  point outside of $H$. Let $t_0$ be different from either of $t_i$
($i=1,...,N$). Let
$\gamma_i$ ($i=1,...,N$) be a good collection, in the sense described
above (Def.(2.1)), of  paths in $\Omega$ based in $t_0$ and forming a
basis of
$\pi_1 (\Omega- \bigcup_i t_i,t_0)$ and let
$\Gamma _i \in Aut (\pi_n ({\bf C}^n_t-V \cap {\bf C}^n_t))$   be the
monodromy automorphism corresponding to $\gamma_i$. Let $\sigma_i:
\pi_{n-1} ({\bf C}^n_{t_i}-V \cap {\bf C}^n_{t_i}) \rightarrow
\pi_n ({\bf C}^n_{t_0} - V \cap {\bf C}^n_{t_0}) ^{\Gamma_i}$ be the
 degeneration operator of the homotopy group of a special
 element of the pencil
into the correponding quotient of covariants constructed above. Then
 $$\pi_n ({\bf C}^{n+1}-V \cap {\bf C}^{n+1})= \pi_ n ({\bf C}^n-V \cap
 {\bf C}^n) / (Im (\Gamma_1-I),Im \sigma_1,...,
 Im (\Gamma_N -I), \sigma_N )\eqno (2.5)$$
\bigskip
\par {\bf Proof.} Let $P: {\bf C}^{n+1} \rightarrow \Omega$ be the
projection defined by the pencil ${\bf C}_t$.
 Let $T({\bf C}_{t_i}^n)$ be the intersection of the
tubular neighbourhood of $L_i$ in ${\bf CP}^{n+1}$ with the finite part
${\bf C}^{n+1}$ which can be taken as $P^{-1} (\Delta _i)$ where
$\Delta_i \subset \Omega$ is a small disk about $t_i$ ($i=1,...,N)$.
 Each loop $\gamma_i$ is isotopic to the loop having the standard form
$\delta^{-1} _i \circ \partial \Delta_i \circ \delta_i$  ($\delta _i $,
as above, is  a system  of  paths in $\Omega$ connecting $t_0$ to
$\partial \Delta _i$ and non-intersecting outside of $t_0$).
 We shall assume from now on that $\gamma_i$'s
have such form. The restriction $P_V$ of $P$ on ${\bf C}^{n+1}-V \cap {\bf
C}^{n+1}$
 defines over $\Omega - \bigcup _i (\gamma_i \cup \Delta_i)$ a locally
trivial fibration
and therefore ${\bf C}^{n+1}- V \cap {\bf C}^{n+1}$ is homotopy
equivalent to $P_V ^{-1} (\bigcup_i (\gamma_i \cup \Delta_i))$.
The latter space is homotopy
equivalent to $$\bigcup_{T({\bf C}_{t_0}^n)- V \cap T({\bf C}_{t_0}^n)}
 T({\bf C}_{t_i}^n)-V \cap T({\bf C}_{t_i}^n) (i=1,...,N)\eqno (2.6)$$
  with the embedding of the
common part of the spaces in the union in each of them depending on the
trivialization of $P_V^{-1} (\delta _i)$ over $\delta_i$. We are
going to calculate the homology of the infinite cyclic cover of
(2.6) by repeated use of  the Mayer-Vietoris sequences.
 First we claim that if $t_{0,i}$ denotes the end point
 of the path $\delta_i$ and $\Gamma^{\prime}_i$ is  the authomorphism of
 $\pi _n ( {\bf C}^n_{t_{0,i}}- {\bf C}^n_{t_{0,i}} \cap V )$
 induced by the monodromy corresponding to the loop $\partial \Delta_i$
 then:
$$ H_n (\widetilde {(T({\bf C}_{t_i}^n -V \cap T({\bf C}_{t_i}^n)},{\bf
Z})=\pi_n ({\bf C}_{t_{0,i}}^n -V \cap
{\bf C}_{t_{0,i}}^n)/(Im(\Gamma^{\prime}_i-I),Im \sigma^{\prime}_i)
\eqno (2.7)$$
 for any $i$ $(i=1,...,N)$ where $\widetilde X$
 denotes the universal cyclic cover of a space $X$
. To verify (2.7) let us consider a small polydisk $B_i \subset T({\bf
C}^n_{t_i}$ for which the projection $P$ induces a split
 $B_i=\Delta^\prime _i \times D^n_i$ as a product of
a 2-disk $\Delta^\prime _i \subset \Delta_i$, $t_i \in \Delta
^\prime_i$ and
$n$-disk $D^n_i$ in ${\bf C}_{t_i}^n$ such that $D^n_i$ contains
the singular point of $V \cap {\bf C}_{t_i}^n$. One has a natural
retraction $B_i - B_i \cap V$ onto $\partial B_i -\partial B_i \cap V$
which shows that $T({\bf C}_{t_i}^n)-T({\bf C}_{t_i}^n) \cap V$
is homotopy equivalent to $T({\bf C}_{t_i}^n)-T({\bf C}_{t_i}^n)
 \cap V-B$. Let us decompose the latter as :
$$P_V^{-1} (\Delta_i -\Delta^{\prime}_i) \cup (P_V^{-1}
 (\Delta^{\prime}_i)-B_i) \eqno (2.8)$$
 The first component in this union, which we shall call
$\Theta_1$, is a locally trivial
fibration over homotopy circle $\Delta_i- \Delta^\prime _i$ with fibre
 ${\bf C}_{t_0}^n-V \cap {\bf C}_{t_0}^n$. The second, which we denote
$\Theta _2$, is fibred over
2-disk $\Delta^\prime _i$ with the fibre ${\bf C}_{t_i}^n-{\bf C}_{t_i}^n
 \cap V$ and hence is homotopy equivalent to this fibre.
  The intersection $\Theta_0$ of
 two pieces in (2.8) which is the  preimage of the boundary circle
  $\partial \Delta ^ \prime_i$ forms a part of
a fibration over a disk and hence is homotopy equivalent to
${\bf C}_{t_i}^n-{\bf C}_{t_i}^n \cap V$. This decomposition defines the
decomposition of the infinite cyclic covers:
$$\widetilde {T({\bf C}^n_{t_i})-T({\bf C}^n_{t_i} \cap V)}=
 \tilde \Theta_1 \bigcup _{\tilde \Theta_0} \tilde \Theta_2
  \eqno (2.9)$$
 The split of $\Omega_0$ as a direct product ${\bf C}^n_{t_i}
 -{\bf C}^n_{t_i} \cap V \times S^1$
 implies the splitting of the infinite cyclic cover of $\tilde
 \Omega_0$ as $\widetilde {({\bf C}^n_{t_i}-{\bf C}^n_{t_i} \cap V)}
 \times S^1$. Therefore $$H_j(\tilde \Omega_0,{\bf Z}))= H_{j-1}
 (\widetilde {{\bf C}^n_{t_i}-{\bf C}^n_{t_i} \cap V},{\bf Z})
 \oplus H_j(\widetilde {{\bf C}^n_{t_i}-{\bf C}^n_{t_i} \cap V},
 {\bf Z}) (j \in {\bf Z})\eqno (2.10)$$
 Next let us consider the Mayer-Vietoris homology sequence corresponding
 to decomposition (2.6): $$\rightarrow H_n(\tilde \Omega_0,{\bf Z})
\rightarrow H_n (\tilde \Omega_1,{\bf Z})
 \oplus H_n(\tilde \Omega_2,{\bf Z})
  \rightarrow H_n (\widetilde {T({\bf C}^n
  _{t_i})-T({\bf C}^n_{t_i}) \cap V},{\bf Z})\rightarrow$$
$$ \rightarrow  H_{n-1} (\tilde \Omega_0) \rightarrow H_{n-1}
(\tilde \Omega_1)  \oplus H_{n-1} (\tilde \Omega_2) \eqno (2.11)$$
 Th group $H_n (\tilde \Omega_1,{\bf Z})$ can be identified as above
with $\pi _n ({\bf C}^n_{t_{0,i}}-{\bf C}^n_{t_{0,i}} \cap V)
 /(\Gamma^{\prime}_i-I)$ and the left
 homomorphism in (2.11) takes the second summand in (2.10) in the case
 $j=n$ isomorphically to $H_n (\tilde \Omega_2,{\bf Z})$. Hence the
$Coker$ of the left homomorphism  in (2.11) coincides with $\pi _n
  ({\bf C}^n_{t_{0,i}}-{\bf C}^n_{t_{0,i}} \cap V)/(\Gamma^{\prime}_i
 -I,Im \sigma^{\prime}_i)$. Moreover the same argument shows that
 the homomorphism  of $H_{n-1} (\tilde \Omega_0,{\bf Z})$  in (2.11) is
 an injection because $H_{n-2} (\widetilde {{\bf C}^n_{t_{0,i}}-
 {\bf C}^n_{t_{0,i}}},{\bf Z})$ is isomorphic to $\pi _{n-2}$ of
the same space  and therefore is trivial (cf. Lemma (1.5)) which implies
(2.7).
\par  Next we shall calculate the homology of the infinite cyclic cover
 in (2.6) inductively using the Mayer-Vietoris sequence
corresponding
 to this decomposition. Because
 $\pi_{i} ({\bf C}^n_{t_0}-{\bf C}^n_{t_0} \cap V)=0$ for
 $2 \le i \le{n-1}$ (cf. lemma
 (1.5)), the terms in the Mayer-Vietoris sequence below dimension $n$
 vanish. The cokernel of the map:
 $$H_n (\widetilde {T({\bf C}^n_{t_0})-T({\bf C}^n_{t_0}) \cap V},
  {\bf Z})\rightarrow$$
 $$\rightarrow
  H_n (\widetilde {{\bf C}^n_{t_{0,i}}-{\bf C}^n_{t_{0,i}}
 \cap V},{\bf Z})/(Im (\Gamma^{\prime}_i-I),Im \sigma^{\prime}_i)
 \oplus H_n (\widetilde {{\bf C}^n_{t_{0,j}}-{\bf C}^n_{t_{0,j}}
 \cap V},{\bf Z})/(Im (\Gamma^{\prime}_j-I),Im \sigma^{\prime}_j)
 \eqno (2.12)$$
 by the linear algebra is isomorphic to $H_n (\widetilde
 {{\bf C}^n_{t_0}-{\bf C}^n_{t_0} \cap V},{\bf Z})/ Im (\Gamma_i-I),
 Im \sigma_i, Im (\Gamma_j -I),Im \sigma_j)$.
 Now the theorem follows.
  \bigskip
\bigskip

\ssectitle{3.\  A vanishing theorem.}

\bigskip
\par In this section we give a necessary condition for the vanishing of
 $\pi_n ({\bf C}^{n+1}-V)$. This is useful in applications of theorem 2.4
 because it allows one to dispose of
  degeneration operator in some cases. We
 give a numerical consequence in the case of curves. The key part is the
following counterpart of the commutativity of the fundamental group in
the case of curves.
\bigskip
\par {\bf Theorem 3.1.} Let $V$ be a hypersurface in ${\bf CP}^{n+1}$
which has only isolated singularities including singularities at
 infinity, and let $H$ be the hyperplane at infinity.
  Suppose that there is a non
singular variety $W$ and a map $\phi: W \rightarrow {\bf CP}^{n+1}$ such
that the union of the proper preimage $V'$ of $V$,the proper
preimage $H'$ of $H$, and the exceptional locus $E$ of $\phi$, form
 a divisor with normal crossings on $W$. Assume that $V'$ is an
{\it ample}  divisor on $W$.
 Then the action of $\pi _1 ({\bf CP}^{n+1}-(V \cup H)={\bf Z}$ on
$\pi _n ({\bf CP}^{n+1}- (V \cup H))$ is trivial.
\bigskip
\par {\bf Proof.} Let $U_{V'}$ be a tubular neighbourhood of $V'$ in
 $W$. First note that the map $\pi _i ((U_{V'}-V')- (E \cup H') \cap
 (U_{V'}-V')) \rightarrow \pi_i (W -  (E \cup V' \cup H'))$
 induced by inclusion is an isomorphism for $i \le n-1$ and
 is surjective for $i=n$. This follows immediately from
  the assumption of ampleness of $V'$ and the Lefschetz theorem
  for open varieties (cf. [H]. In the setting of this work one
  applies the theorem 2 from this paper  to $W$ embedded into ${\bf
CP}^N$ using a  multiple of the line bundle corresponding to $V'$ and
taking $V'$  as the hyperplane section at infinity).
 Because $W- (V' \cup H' \cup E)=
 {\bf CP}^{n+1} -(V \cup H)$  we obtain $\psi _i: \pi_i ((U_{V'}-V')
  -(E \cup H') \cap (U_{V'}-V')) \rightarrow \pi_i ({\bf CP}^{n+1}
  - (V \cup H)$ which is an isomorphism for $i \le n-1$ and surjective
 for $i=n$.
 \par Let $\alpha$ be the boundary of the normal to $V'$ 2-disk in
 $U_{V'}$ at a point of $V'$ outside of $E \cup H'$. The action of
 $\alpha$ considered as an element of $\pi_1 ((U_{V'}-V')
 -(E \cup H') \cap (U_{V'}-V'))$ on $\pi_n ((U_{V'}-V')-
 (E \cup H') \cap (U_{V'}-V'))$ is trivial. Indeed $(U_{V'}-V')
 -(U_{V'}-V') \cap (E \cup H')$ is a (trivial) circle bundle
 over $V'-(V' \cap (E \cup H')$ because $V' \cup H' \cup E$ is assumed
 to be a divisor with normal crossings. Moreover the projection map
 $\lambda$ induces the isomorphism $\lambda_*: \pi_i ((U_{V'}-V')
 - (U_{V'}-V') \cap (E \cup H')) \rightarrow \pi _i (V'-V' \cap H'
 \cap E)$ for $i \ge 2$. If $\gamma \in \pi_i ((U_(V')-V')- (U_(V')
 -V') \cap (E \cup H'))$ then $\lambda _* (\alpha \cdot \gamma)=
  \lambda _*(\alpha) \cdot \lambda _* (\gamma) = \lambda_* (\gamma)$
 i.e. $\alpha \cdot \gamma = \gamma$ and our claim follows. This
also concludes the proof of the theorem because
  $\psi_n$ is surjective and because $\alpha$ is taken by $\psi_1$
into the
 generator of $\pi _1 ({\bf CP}^{n+1}- (V \cup H))$.
\bigskip
\par {\bf Theorem 3.2.} Let $V$ be a hypersurface in ${\bf CP}^{n+1}$
 which satisfies all conditions of the theorem 3.1. Let us assume also
 that $H^{n+1} (V, V \cap H,{\bf Z})=0$. Then $\pi _n ({\bf CP}^{n+1}
 - (V \cup H))$ vanishes.
\bigskip
\par {\bf Proof}. Let us consider the Leray spectral sequence
$$E_{p,q}^2 =H_p ({\bf Z},
 H_q (\widetilde {({\bf CP}^{n+1}-(V \cup H)},{\bf Z}))
 \Rightarrow H_{p+q} ({\bf CP}^{n+1}- (V \cup H),{\bf Z}) \eqno (3.1)$$
  asssociated
 with the classifying map of ${\bf CP}^{n+1}-(V \cup H)$ into
  $S^1=K({\bf Z},1)$ corresponding to the generator of $H^1({\bf CP}
 ^{n+1}- (V \cup H),{\bf Z})={\bf Z})$  (the homotopy fibre of this map
 is the universal cyclic cover $\widetilde {({\bf CP}^{n+1}-(V \cup
H))}$. It  implies that $H_0 ({\bf Z},H_n {(\widetilde ({\bf CP}^{n+1}
 -(V \cup H))},{\bf Z})$ which isomorphic to
the covariants $H_n (\widetilde {({\bf CP}^{n+1}-(V \cup H))}
  ,{\bf Z})^{\pi_1 ({\bf CP}^{n+1}- (V \cup H))}=H_n ({\bf CP}^{n+1}
 -(V \cup H),{\bf Z})$. The group in the left hand side is isomorphic to
  the covariants $\pi_n ({\bf CP}^{n+1}-(V \cup H))^{\pi _1 ({\bf CP}
 ^{n+1}- (V \cup H))}$. Hence the result follows from the above theorem,
 and the vanishing of $H_n({\bf CP}^{n+1}- (V \cup H),{\bf Z})$ which
 is a consequence of the lemma 1.6 and the assumption we made on the
cohomology of $(V,V \cap H)$.
\bigskip {\bf Corollary 3.3.} Let $V$ be a hypersurface in ${\bf
CP}^{n+1}$ which satisfies the condition of the theorem 3.1. If $V$ and
$V \cap H$ are ${\bf Q}$-manifolds then $\pi_n({\bf CP}^{n+1}-(V \cap
H))=0$.
\bigskip {\bf Proof.} This follows from 3.2 and 1.7.
\bigskip {\bf Remark 3.3} In the case $n=1$ the argument given in the
 proof of the theorem 3.1 can be strengthened to show that $(V \cdot V$
 implies that the fundamental
group of the complement is abelian (cf. [Ab],[N]).
 Recall that, for example for a curve $C$ of degree $d$ which has
$\delta$ nodes and $\kappa$ cusps, this implies the commutativity of
 $\pi_1 ({\bf CP}^{n+1}-C)$ provided $4 \delta + 6 \kappa <d^2$.  For
the application of the technique of section 2 the following result is
useful.
\bigskip {\bf Lemma 3.4} Let $C$ be a curve in ${\bf CP}^2$ which has
only one singular point which is unibranched and has one characteristic
pair $(m,k)$ ($k \le m$). If $d^2 >m \cdot k$ then $\pi_1 ({\bf
CP}^2-C)$ is abelian. In particular  the Alexander module of $C$ is
trivial.
\bigskip
\par {\bf Proof.} The greatest common divisor of $m$ and $k$ which
is equal to the number of branches of the singularity is equal to 1. The
resolution of singularities of plane curves can be described in terms of
the Euclidian algorithm for finding this greatest common divisor
(cf.[BK]). Let $m=a_1k+r_1,k=a_2 \cdot r_1+r_2,...r_{s-1}=a_{s+1} \cdot
r_s +1$ be the steps of the Euclidain algorithm. Then the sequence of
blow ups which results in the embedded resolution produces the
following. Each of the first $a_1$ blow ups  gives an exceptional
curve
with multiplicitity $k$, the intersection index of which with the proper
 preimage of the curve in quesiton is equal to $k$. This results
in dropping of
the self-intersection index of the proper preimage by $k^2$.
Subsequent blow ups have a similar effect with the blow ups
corresponding to the last step in Euclidian algorithm resulting in a
non-singular proper preimage with the tangency order with the exceptional
 curve equal to $r_s$. Additional $r_s$ blow ups result in a  proper
preimage, the union of which with the exceptional locus, is  a divisor
 with normal crossings. The self intersection of the proper preimage is
equal to $d^2-a_1 \cdot k^2-a_1 \cdot r_1^2-...-a_{s+1} \cdot r_s^2
-r_s=d^2-(m-r_1) \cdot k-...-(r_{s+1}-1) \cdot r_s -r_s=d^2-m \cdot k
>0$. Hence our claim follows from the Nori's theorem [N].
\bigskip
\bigskip

\ssectitle{4.\ The divisibility  theorems}.

\bigskip
\par In this section we prove two theorems relating the order of
 the homotopy group of the complement to the hypersurface $V$
 to the type of the singlarities $V$ including the singularities
of $V$ at infinity. We will discuss the relation of these
results to the divisibility theorem for Alexander polynomials in [L].
 We asssume as above that $V$ has only isolated singularities.
\par First recall that
if $c \in V$ is a singular point of $V$ then one can associate with it
the  characteristic polynomial $P_c$
of the monodromy operator in the Milnor
 fibration of the singularity $c$. By a certain abuse of language
we will call this polynomial the polynomial of the singularity $c$.
 The cyclic cover $U_c$ of $B_c- B_c \cap V$ corresponding to
  the kernel of the homomorphism
  $lk_c: \pi_1 (B_c-B_c \cap V) \rightarrow
 H_1(B_c-B_c \cap V,{\bf Z})={\bf Z}$ (it is given by the
 evaluation of the linking number with $V \cap B_c$) is
  homotopy equivalent to the Milnor fibre
 and the characteristic polynomial of the automorphism of
 $H_n(U_c,{\bf Q})$ induced by the deck transformation
coincides with the polynomial of the singularity $c$ as a consequence
 of existence of fibration of $B_c-B_c \cap V$.
 On the other hand
 if $c \in Sing_{\infty}$
 and $B_c$ is a small ball in a certain Riemannian metric in
 ${\bf CP}^{n+1}$ about $c$, then $H_1(B_c-(V \cup H) \cap B_c,{\bf Z}
 )= H_2(B_c,B_c-(V \cup H),{\bf Z})=H^{2n}(T(V \cup H), S_c \cap T(V \cup
 H)$ where $S_c$ is the boundary of $B_c$ and $T()$ denotes the regular
 neighbourhood. The latter is isomorphic to $H^{2n}(V \cap B_c, \partial
 (V \cap B_c),{\bf Z}) \oplus H^{2n} (H \cap B_c,\partial (H \cap B_c),
 {\bf Z})={\bf Z} \oplus {\bf Z}$. The map $lk_c :
 \pi _1 (B_c- (V \cup H)
 \cap B_c) \rightarrow {\bf Z}$ corresponding to projection of
  $H_1 (B_c- B_c \cap (V \cup H)$ onto $H^{2n} ( B_c \cap V ,\partial
(B_c \cap V), {\bf Z})={\bf Z}$ geometrically is the linking number with
 $V$. Let  $U_ c$  be the infinite cyclic cover of $B_c-(V \cup H)
 \cap B_c$ corresponding to the kernel of the homomorphism $lk_c$.
\bigskip
 \par {\bf Definition 4.1.} The order $P_c(t)$
 of $H_n(U_c,{\bf Q})$ considered
 as the module over ${\bf Q} [t,t^{-1}]$ via the action induced by
  the deck transformations on $U_c$ is called the polynomial
  of the singular point $c$.
\bigskip
\par {\bf Remark 4.2.} $U_c$ has a homotopy type of a finite dimensional
 complex and in particular $H_i(U_c,{\bf Q})$ is a torsion
 ${\bf Q}[t,t^{-1}]$-module for any $i$. So the above polynomial
 $P_c(t)$ is non zero. This follows from the following realization of
  the infinite cyclic cover $U_c$.
 Let $\phi_c$ be a holomorphic function in
 a neighbourhood $N_c$ of $c$ in ${\bf CP}^{n+1}$ such that $\phi_c=0$
 coincides with $V$ in this neighbourhood. Let $V_c(s)$ be given by
 the equation $\phi_c=s$ in $N_c$. Then $U_c$ is homotopy equivalent to
   $V_c(s)-H \cap V_c(s)$ for $s$ sufficiently close to zero. Indeed
 the union of hypersurfaces $\phi_c=s$ ($s \le \epsilon$ and $N_c$ is
 sufficiently small) is homeomorphic to the ball ([Mi]) and the function
 $\phi_c$ provides the locally trivial fibration of this ball
  over a punctured disk. Because the
 singularity of $V \cap H$ is isolated $H$ will be transversal to
 all hypersurfaces $\phi_c=s$ ($\epsilon$ sufficiently small). Therefore
 $\phi_c(s)$ also provides the locally trivial fibration of the
 complement in this ball to $V \cap H$. Hence
 $V_c(s)-V_c(s) \cap H$ is homotopy equivalant to $U_c$.
\bigskip
 \par {\bf Theorem 4.3.} If $V$ and $V \cap H$ have only isolated
 singularities then the order $P_V$
 of $\pi_n ({\bf CP}^{n+1}- (V \cup H)) \otimes {\bf Q}$
 as the module over ${\bf Q}[t,t^{-1}]$ divides the product
  $\Pi_c P_c \cdot (t-1)^{\kappa}$ (for some $\kappa \in \bf Z$
  of the polynomials $P_c$ of all
 singularities of $V$ including those in
  ${\bf C}^{n+1}={\bf CP}^{n+1}-H$ as well as those at  infinity.
 One can drop the factor $(t-1)^{\kappa}$ if one of the following
 conditions takes place:
\par a) $H^{n+1}(V,H \cap V,{\bf Q})=0$.
\par b) $\pi_n ({\bf CP}^{n+1}-(V \cap H) \otimes \bf Q$ is
semisimple ${\bf Q} [t,t^{-1}]$-module.
\bigskip
\par {\bf Remark 4.4.} One of the consequences
 of the theorem is that the order
 of $\pi _n ({\bf CP}^{n+1}-(V \cup H))$ is not zero. Hence this
 ${\bf Q}[t,t^{-1}]$-module is a torsion module.
\bigskip
\par {\bf Proof}. Let $T(V)$ be a regular neighbourhood of
 $V$ in ${\bf CP}^{n+1}$. First let us observe that
  $\pi_n (T(V)-((T(V) \cap (H \cup V))$ surjects onto
 $\pi _n({\bf CP}^{n+1} - (V \cap H))$. Indeed $T(V)$ contains
 a generic hypesurface $W$ of the same degree as $V$
 which is transversal to $V \cup H$.
 According to the Lefschetz theorem  $\pi _i(W- W \cap (V \cup H))$ maps
 by the map $j_*$ induced by inclusion isomorphically
 to $\pi_i({\bf CP}^{n+1}- (V \cup H))$ for $1 \le i \le {n-1}$ and
 surjects for $i=n$. Now our claim follows from the fact that $j_*$
 can be factored as
 $$\pi_n(W-W \cap (V \cup H)) \rightarrow
\pi_n(T(V)-T(V) \cap (V \cup H)) \rightarrow \pi_n ({\bf CP}^{n+1}
- (V \cup H)). \eqno (4.1) $$.
\par  Next let us consider a   collection of non intersecting
 balls  $B_c$ about the singular points of
  $V$, $c \in Sing (V) \cup  Sing _\infty (V)$. For a sufficiently
small
 regular neighbourhood  $T(V)$ of $V$ the complement to the union of
 $B_c$ ($c \in Sing (V) \cup Sing _{\infty} (V)$):
  $$B_0= T(V) -(T(V) \cap (V \cup H))-
 \bigcup_c  B_c \cap (T(V)-T(V) \cap (V \cup H)). \eqno (4.2) $$ fiberes
over
 $V-(H \cap V  \bigcup _c (B_c \cap V))$ with the circle as the fibre.
 Such a fibre maps onto the generator of
 $H_1 ({\bf C}^{n+1}- V,{\bf Z})$. Hence
 one has the surjection $lk_ T :H_1(T(V)-T(V) \cap (V \cup H),{\bf
Z})
 \rightarrow H_1({\bf C}^{n+1}-V,{\bf Z})={\bf Z}$ ( the linking
number, cf. Lemma 1.6). The kernel of $lk_T$ defines the infinite
 cyclic cover $U_T$ of $T(V)-T(V) \cap (V \cup H)$.
For any $c$ the
 map $lk_c$ coincides with the composition of the map of the
 fundamental groups induced by
embedding $B_c - B_c \cap V \rightarrow T(V)-T(V) \cap (V \cup H)$
 and the map $lk_T$. A similar factorization takes place for
 the linking number homomorphism of $\pi_1 (B_0)$, which defines
 the infinite cyclic cover $U_0$ of $B_0$.
 We obtain that $$U_T =U_0 \bigcup_c U_c. \eqno (4.3)$$
 We claim that $U_T$ is $n-1$-connected and in particular
 $\pi_n (T(V)- T(V) \cap (V \cup H) =H_n(U_T,{\bf Z})$.
 Indeed the fundamental group of $U_T$ can be obtained from the
 fundamental group of $U_0$ by induction corresponding to adding
 $U_c$ one by one using the van Kampen theorem (on $\pi_1$ of the union)
 . Each time the
 fundamental group is replaced by the quotient by the image
of the fundamental group of the link of the corresponding singularity
 (in the case $c \in Sing_\infty (V)$ one rather should take the
 quotient by the image of the fundamental group of the complement to
 the intersection of the link of the singularity with the hyperplane
 at infinity inside this link). But the fundamental group of the affine
 portion of
 the smoothing of $V$ which is transversal to $H$ is calculated by
 the same procedure. Because the latter is simply connected we obtain
 that the fundamental group of $U_T$ is trivial. On the other hand we
 have the following Mayer-Vietoris sequence:
 $$\oplus _c H_i(U_0 \cap U_c) \rightarrow \oplus _c  H_i (U_c) \oplus
 H_i (U_0) \rightarrow H_ i (U_T) \rightarrow. \eqno (4.4)$$
  For $2 \le i \le n-1$ we have $H_i (U_c)=H_i (V_c(s))=0$.
 (Here, as in remark 4.2,  $V_c (s)$
  is a smoothing of the singularity $c$).
 This follows from the standard connectiviy of the Milnor
 fibre for finite singularities and for singularities at infinity
 follows from the latter  and the exact sequence of the pair and
 $H_i (V_c(s)- V_c(s) \cap H,V_c(s))=
  H_{i-2} (V_c (s) \cap H)=0$ for $2 \ne i \le n$ (the first
 isomorphism is a consequence of excision and the Poincare duality).
 On the other hand if $V'$ is a hypersurface in a pencil of
 hypersurfaces which contains $V$ and has $V \cap H$ as the
 base locus then $H_i (V', V' \cap H, \cup_c V_c (s))=
 H_i (V',V' \cap H)=0$
 ($c \in Sing (V)
 \cap Sing_ \infty (V)$) for $i \le n-1$ by Lefschetz theorem.
 Now excision of the union of small balls about all singular points
of $V$ shows that $H_i (U_0 \cap (\cup _c U_c))=
H_i (\cup_c \partial V_c (s), \partial V_c (s) \cap H)
  \rightarrow H_i (U_0)= H_i (V'- \cup V_c(s)
,V' - \cup V_c(s) \cap H) $ is
surjective for $i \le n-1$ and is isomorphism for  $i \le n-2$.
Because $H_i (U_0 \cap U_c)=H_i (\partial V_c (s))=0$ for $0 < i <
{n-1} $ we obtain that $H_i (U_T)=0$ for $0 <i <n$.
\par For $i=n$ the sequence (4.4), being equivariant with respect to the
action of ${\bf Z}$ by the deck transformations implies that
if  $ Q=ord (Ker \oplus_c H_{n-1} (U_0 \cap U_c) \rightarrow \oplus_c
H_{n-1} (U_c) \oplus H_{n-1} (U_0))= ord
(Ker \oplus_c (H_{n-1} (U_0 \cap U_c)
 \rightarrow H_{n-1} (U_0))$ and $R=ord (Ker \oplus_c H_n
 (U_0 \cap U_c) \rightarrow \oplus_c H_n (U_c)
\oplus H_n (U_0))$ then $ord(H_n (U_T)= Q \cdot R$. The orders of
$H_{n-1} (U_0 \cap U_c)$ and $H_n (U_0)$ are powers of $t-1$
because those are the cyclic covers of the (trivial) circle bundles.
 The order of $\oplus H_n(U_c)$ is equal to $\Pi_c P_c$ and hence $R$
 divides this product multiplied by a power of $t-1$.
 Therefore  the order of $H_n (T(V)-T(V) \cap (V \cup H)$
  divides $\Pi_c P_c$ multiplied by $(t-1)^ {\kappa}$ for some $\kappa$.
It follows from (4.1) that the same is true for order $P_V(t)$ of $\pi
_n ({\bf CP }^{n+1} -(V \cup H)$.
 To conclude the proof we need to show that the order of zero of $P_V
(t)$ at 1 does not exceed the sum of the orders of the zero at 1 of
 $P_c (t)$.
\par If one assumes a) above then according to the lemma  1.12 and 1.6
 $P_V (1) \ne 0$ and the theorem follows. Moreover in the case b)
the order of the zero of $P_V (t)$ at 1 is equal to the
rank of $H^{n+1} (V,H \cap V,{\bf Q})$ as follows from the sequence
 (1.3).
 The order of the zero
 of $P_c (t)$ at 1 is greater or equal than the rank $H_{n-1} (L_c,{\bf
 Q})$ where $L_c$ is the link of the singularity $c$ (cf. [Mi]).
 Hence in the case b) the theorem follows from the inequality:
 $$rk H^{n+1} (V,V \cap H,{\bf Q}) \le \Sigma _c  rk H_{n-1}
 (L_c ,{\bf Q}) \eqno (4.5)$$
 To show this, note that $H^{n+1}(V,V\cap H, {\bf Q})=
 H^{n+1} (V,V \cap H \cup S,{\bf Q})=
 H^{n+1} (V- T(H) \cap V,(\partial T(H) \cap V) \cup T(S),{\bf Q})$
 where $S$ is the collection of the singular points of $V$ outside $H$,
 $T(S)$ is a small regular neighbourhood of this finite set in $V$
 and $T(H)$ is the regular neighbourhood of $V \cap H$ in $V$. The
 last cohomology group is
  dual to $H_{n-1} (V-V \cap H-S,{\bf Q})$ (use excision of $T(S)$).
  (4.5) will follow from the exact sequence of the pair and the
vanishing of $H_{n-1} (V-S-V \cap H, \cup L_c, {\bf Q})$. This
 group is isomorphic to $H_{n-1} (\tilde V-H \cap V, \cup M_c, {\bf Q})$
 where $M_c$ is the Milnor fibre of the singularity $c$ and $\tilde V$
 is a generic hypersurface in the pencil of hypersurfaces containing
$V$ and having $V \cap H$ as the base locus. The vanishing of
 the last group is a consequence of the $n-1$ connectedness of the
 Milnor fibres $M_c$ and the vanishing of $H_{n-1} (\tilde V - V \cap H)$
  follows from the arguments just used and the exact sequence of the
pair $ (\tilde V, \tilde V -H \cap V)$.
\bigskip
\par {\bf  Theorem 4.5.}
 Let $V$ be a hypersurface in ${\bf CP}^{n+1}$ having
only isolated singularities including infinity. Let $H$ be the
hyperplane at infinity. Let $S_\infty$ be a sphere of a sufficiently
large radius in ${\bf C}^{n+1} ={\bf CP}^{n+1}-H$ (or equivalently the
boundary of a suffuciently small tubular neighbourhood of $H$ in ${\bf
CP}^{n+1}$). Let $U\infty $ be the infinite cyclic cover of $S_\infty
-V \cap S_\infty$  corresponding to the kernel of the
 homomorphism: $\pi _1(S_\infty -V \cap S_\infty)
\rightarrow {\bf Z}$ given by the linking number(cf. remark below). Let
 $P_\infty$ be the order of $H_n(U_\infty, {\bf Q})$ as the ${\bf Q}
[t,t^{-1}]$-module. Then $P_V$ divides $P_\infty$.
\bigskip
\par {\bf Remark 4.6}. $V \cap S_\infty$ is a connected manifold if
$n \ge 2$. If $n=1$, the number of connected
components of $V \cap S_\infty$ is "the number of places of
the curve at infinity. By Alexander duality if $H_1 (S_ \infty -V \cap
S_ \infty,{\bf Z})=H^{2n-1} (V \cap S_ \infty,{\bf Z})
={\bf Z}$ if $n \ge 2$
and for curves $H_1(S_\infty -S_\infty \cap V,{\bf Z})$ is a free
abelian group of the rank equal to the number of places at infinity.
\bigskip
\par {\bf Remark 4.7.} $H_n (U_\infty ,{\bf Q})$ is a torsion module.
 Indeed let $L_c(V)$ (resp. $L_c(V \cap H)$) be the link of the
 singularity  $c$ of $V$ (resp. $V \cap H$) in ${\bf CP}^{n+1}$
 (resp. $H$). Let $B_c$ be a polydisk in ${\bf CP}^{n+1}$ of the form
 $D^{2n}_c \times D^2_c$ about  $c$ such that the part of its boundary
 $S^1_c \times D^{2n}_c \subset \partial T (H)=S_\infty$. Then
 $$ S_\infty-S_\infty \cap V=\partial T(H)-\partial T(H) \cap V
 =\bigcup_c(S^1_c \times D^{2n}_c -V) \cup U \eqno (4.6)$$
 where $U$ fibres over $H-H \cap V$ with a circle as a fibre and hence
 the homology of the infinite cyclic cover of $U$ is a torsion
 ${\bf Q}[t,t^{-1}]$-module (of the order which is a power of
 $t^d-1$ where $d$ is the degree of $V$  because the circle
 about $H$ is homologous to $d$ multiplied by the generator
 of $H_1 ({\bf CP}^{n+1}- (V \cup H))$).
  $S^1_c \times D^{2n}_c -V=B_c-(V \cup H) \cap B_c$ and the
 homology groups of the infinite  cyclic cover in question
 are ${\bf Q} [t,t^{-1}]$-torsion modules as follows from
 remark 4.2. Finally the intersection in the union (4.6) of $U$
 with $S^1_c \times D^{2n}_c -V \cap S^1_c \times D^{2n}$
 fibres over $L_c(V \cap H)$  and hence the homology of its
 infinite cyclic cover is a torsion ${\bf Q} [t,t^{-1}]$ module
 as well. Hence the Mayer-Vietoris sequence yields the claim.
\bigskip
\par {\bf Proof of 4.5 }. First note that $S_\infty \cap V$ is homotopy
 equivalent to $T(H)-T(H) \cap (V \cup H)$ where $T(H)$ is the tubular
neighbourhood of $H$ for which $S_\infty$ is the boundary. If $L$ is a
generic hyperplane in ${\bf CP}^{n+1}$, which we will can assume belongs
 to $T(H)$, then once again by Lefschetz theorem  we obtain that the
composition:
$$\pi_n (L-L \cap (V \cup H) \rightarrow \pi_n (T(H)-T(H) \cap (V \cup
H)) \rightarrow \pi _n ({\bf CP}^{n+1}-(V \cup H)) \eqno (4.7)$$
 is surjective. Now the theorem follows from the multiplicativity of the
order in exact sequences.
\bigskip
\par {\bf Corollary 4.8.} If $V$ is a hypersurface transversal to the
hyperplane $H$ at infinity  then $\pi _n({\bf CP}^{n+1}-(V \cup H))
\otimes {\bf Q}$ is a semisimple ${\bf Q} [t,t^{-1}]$-module.
 Any root of the order $P_V$ of the homotopy group is a root of $1$
 of degree $d$.
\bigskip
\par {\bf Proof}. The surjectivity on the right homomorphism in (4.7)
 shows that the claim will follow from the semisimplicity of
$\pi_n (T(H)-T(H) \cap (V \cup H) \otimes {\bf Q})$ as
 ${\bf Q} [t,t^{-1}] $-module. But $T(H)-T(H) \cap (V \cup H)$ is
 homotopy equivalent to the link of the singularity
 ${\cal V}_0: z_1^d+...+z_{n+1}^d=0$  $(d=deg V)$
 in ${\bf C}^{n+1}$ provided $V$ is transversal to $H$. Indeed
  the projective closure of this hypersurface intersects $H$ in a
  non-singular hypersurface which is isotopic to $V \cap H$ and this
isotopy can be extended to a neighbourhood of $H \cap V$. The monodromy
 of ${\cal V}_0$ is semisimple
  (this is the case for any weighted homogeneous
 singularity because, as one can see from the explicit
  description of it (cf. [M]), this monodromy has a finite order).
   The last part in the statement of
 the corollary follows from the Milnor's formula for the characteristic
polynomial of the monodromy of weighted homogeneous polynomials applied
 the singularity ${\cal V}_0$(cf. [M]).
\bigskip
\par {\bf Corollary 4.9.} Let $V$ be a hypersurface in ${\bf CP}^{n+1}$
 given by equation $f=0$. Assume that the singularities of $V$ have
  codimension $k$ in $V$. If $V$
 is transversal to the hyperplane $H$ at infinity (as stratified space)
 then
 $\pi_k ({\bf CP}^{n+1} - (V \cup H)) \otimes {\bf Q})$ is
 isomorphic as ${\bf Q}[t,t^{-1}]$-module to $H_k(M_f,{\bf Q})$
  where $M_f$ is the Milnor fibre
of the singularity  at the origin in ${\bf C}^{n+2}$ of hypersurface
$f(x_0,...,x_{n+1})=0$  with the usual module structure given by
 the monodromy operator.
\bigskip
\par {\bf Proof.} We can assume that  the singularities of $V$ are
 isolated because the general case can be reduced to this
 using Lefschetz theorems (cf. section 1).
First notice that $\pi_n ({\bf CP}^{n+1}-V)$
 is isomorphic to  $H_n (\widetilde {({\bf CP}^{n+1}-V)_d},{\bf Z})$
 where $\widetilde {({\bf CP}^{n+1}-V)_d)}$ is the $d=degV$-fold
 cyclic covering of ${\bf CP}^{n+1}-V$ because $\pi _1$ of the
latter is ${\bf Z} /d {\bf Z}$. This $d$-fold covering
 is analytically equivalent to the affine hypersurface $f=1$
 which is diffeomorpfic to the Milnor fiber $M_f$. The deck
 transformation in this model of $\widetilde {({\bf CP}^{n+1}-V)_d}$
 corresponds to the transformation induced by multiplication of
each coordinate of ${\bf C}^{n+2}$ by  a primitive root of unity
 of degree $d$. It is well known that this is also a description of
 the monodromy of a weighted homogeneous polynomial (cf. [M]). Finally
 according to lemma 1.13 $\pi_n ({\bf CP}^{n+1}-V ) \otimes {\bf Q} =
 \pi _n ({\bf CP}^{n+1}- (V \cup H)) \otimes {\bf Q} / (t^d -1) \cdot
 \pi_n ({\bf CP}^{n+1}- (V \cup H)) \otimes {\bf Q}$ which is
 isomorphic  to $\pi _i ({\bf CP}^{n+1}- (V \cup H)) \otimes {\bf Q}$
 because of the results which are contained in the corollary 4.8.
\bigskip
\par {\bf Remark 4.10} The corollary 4.9
 is an extension to high dimensions of
a   result due to R.Randell [R]  which gives a similar fact for
 Alexander polynomials. In the case of irreducible
 curves the divisibility theorem 4.3 gives a somewhat weaker result than
  the one in [L2] where it is shown that the Alexander polynomial
  divides the product the characteristic polynomials of all
  {\it branches} of the curve in all singular points.
\bigskip
\bigskip

\ssectitle{5.\ Non trivial $\pi_n$}

\bigskip
\par The purpose of this section is to prove two results which allow one
 to construct special classes of hypersurfaces with isolated
singularities for which $\pi_n ({\bf C}^{n+1} -V) \otimes {\bf Q}$ is
non-trivial. We start with the following lemma which may be of
independent interest.
\par {\bf Lemma 5.1.} Let $p(x_1,...,x_{l+1})$ be a polynomial having
a singularity of codimension $k$  at the origin. Let
$x_i=f_i (z_{i,1},..z_{i,{n_i+1}})$, $(i=1,...,{l+1})$ be a collection
 of polynomials all of which we also shall  assume  have at most
  isolated singularities at the origin. Suppose that $n_i \ge k+1$
 for every $i$ for which $f_i$ has singularity at the origin.
   Then the polynomial of $N=\sum _{i=1}^{l+1} (n_i+1)$
variables  $p(f)=p(f_1(z_{1,1},...,z_{1,n_1+1}),...,
 f_{l+1}(z_{{l+1},1}...z_{{l+1},n_{l+1}+1}))$ has the singularity at the
origin of codimension $k$. If $M_p$ and $M_{p(f)}$ denote
the Milnor fibres of the singularities of $p$ and $p(f)$ respectively
 then $H_k (M_p,{\bf Z}) =H_k (M_{p(f)},{\bf Z})$ as ${\bf Z} [t,t^{-1}]$
 modules where the action of $t$ in each case is given by the
 action of the monodromy operator.
\bigskip
\par {\bf Proof.} First note that the use of induction allows one
to reduce
 this lemma to the case when $f_i (z_{i,1}...,z_{i,n_i+1})=z_{i,1}$
 for $i \ge 2$,
 i.e. when the change of variables takes place only in one of $x_i$'s.
Let us select  $\epsilon_1 >0$
and small ball $B_1$ in ${\bf C}^{n_1+1}$ such that the intersection
 of the hypersurface $f(z_{1,1},...,z_{1,n_i+1})=s$  with $B_1$, provided
$0 < \vert s \vert \le \epsilon_1$, is equivalent
to the Milnor fibre of  $f_1$.
 Let us consider a ball $B_0$ centered at the origin ${\bf C}^{l+1}$ of a
 radius less than  $\epsilon_1$. Let $\eta >0$ be
such that for $0 < \vert s \vert < \eta$
the portion of $p(x_1,...,x_{l+1})=s $ which belongs to
$B_0$ is equivalent to the Milnor fibre of $p$. Let $L$ be the
intersection of $M_\eta$ with the  coordinate hyperplane
$x_1=0$ in ${\bf C}^{l+1}$. Finally let us fix a
polydisk $D \subset {\bf C}^{l+n_1+1}$ projection of which
 on subspace $x_2=...x_{l+1}=0$ belongs to $B_1$ and such
 that intersection of it with $p(f)= s$ for $0 < \vert s \vert < \eta$
 is equivalent to the Milnor fiber of $p(f)$. On a part $D'$
 of the polydisk $D$ the formula:
 $$x_1=f_1 (z_{1,1},...,z_{1,n_1+1}),
  x_i=z_{i,1}, ( i=2,...,l+1)  \eqno (5.1)$$
 defines  a holomorphic map $F: D' \rightarrow B$. This map, when
 restricted on
a Milnor fibre $M_{p(f)} \subset D'$ of $p(f)$ which is given by
$p(f)=\eta$, takes $M_{p(f)}$  onto the Milnor fibre of $p$ given by
$p=\eta$. Let $\tilde L$ be the preimage of $L$: $F^{-1} (L)$. The
restriction of $F$ on $M_{p(f)}-\tilde L$ is a locally trivial
fibration: $F: M_{p(f)}-\tilde L \rightarrow M_p-L$. The fibre of this
fibration is equivalent to  the Milnor fibre  $M_{f_1}$
of $f_1 (z_{1,1},...,z_{1,n_1+1})$, This Milnor fibre is
$(n_1-1)$ connected.
 The Leray spectral sequence: $E_{p,q}^2=H_p(M_p-L,H_q(M_{f_1},{\bf Q})
 \Rightarrow H_{p+q} (M_{p(f)}- \tilde L,{\bf Q})$ shows that the
 isomorphism $H_i(M_{p(f)}- \tilde L,{\bf Q})=
  H_i (M_p-L,{\bf Q})$ will take place for
$i$'s which includes $k$ and $k+1$. The following diagram which
compares two exact sequences of pair:
$$\matrix {{H_{k+1}(M_p(f),M_p(f)-\tilde L,{\bf Q})}&\rightarrow&
{H_k(M_{p(f)}-\tilde L,{\bf Q})}&\rightarrow&{H_k(M_{p(f)},{\bf Q})}
&\rightarrow \cr
\downarrow && \downarrow && \downarrow \cr
{H_{k+1}(M_p,M_p-L,{\bf Q})} & \rightarrow &{H_k (M_p-L,{\bf Q})}
 & \rightarrow & {H_k(M_p,{\bf Q})} & \rightarrow \cr }$$
$$\matrix {\rightarrow & {H_k (M_{p(f)},M_{p(f)}- \tilde L,{\bf Q})} &
\rightarrow & {H_{k-1} (M_{p(f)}- \tilde L,{\bf Q})} & \rightarrow \cr
 & \downarrow && \downarrow & \cr
\rightarrow & {H_k(M_p,M_p - L,{\bf Q})} & \rightarrow & {H_{k-1}
(M_p-L,{\bf Q})} & \rightarrow \cr} \eqno (5.2) $$
and the five lemma yields that the isomorphism of our lemma is a
consequence of the isomorphism $$H_i (M_{p(f)}, M_{p(f)}- \tilde L,
{\bf Q})  \rightarrow H_i (M_p,M_p -L,{\bf Q})  for i=k,k+1
\eqno (5.3)$$  Let $T(L)$ (resp. $T(\tilde L)$)
be the regular neighbourhoods of $L$ (resp. $\tilde L$) in $M_p$ (resp.
$M_{p(f)}$) and $\partial T(L)$  (resp. $\partial T(\tilde L)$) be the
boundary of $T(L)$ (resp. $T(\tilde L)$). Then using
 excision one obtains  that
$H_i (M_{p(f)},M_{p(f)}- \tilde L,{\bf Q})=H_i (T(\tilde L),\partial
T(\tilde L),{\bf Q})$ and $H_i (M_p,M_p-L,{\bf Q})=H_i (T(L),\partial
T(L),{\bf Q})$. But $H_i ( T(\tilde L),{\bf Q}) = H_i (T(L),{\bf Q})$.
Indeed $\tilde L$ fibres over $L$ with contractible fibre (i.e. the
 part of $f_1=0 $ inside $B_1$ which is the cone over the link of the
 singularity $f_1$) and $\partial  T(\tilde L)$ is a
fibration over $\partial T(L)$ with $(n_1-1)$-connected fibre. Hence the
 last claimed isomorphism follows from the five lemma.
\bigskip
\par This lemma has the following corollary:
\bigskip
\par {\bf Theorem 5.2.} For an integer $k$ let $g_k (z_0,...,z_{n+1})$
 be a generic form of degree $k$. Let $V_{d_1,...,d_{n+1}}$  be a
hypersurface of degree $D=d_1 \cdot d_2 \cdot ...\cdot d_{n+1}$ given by
the equation:
 $$g_{d_2 \cdot... \cdot  d_{n+1}}^{d_1}+g_{d_1 \cdot d_3 \cdot ...\cdot
d_{n+1}}^{d_2}+....+ g_{d_1 \cdot d_2 \cdot ... \cdot d_n}^{d_{n+1}}=0
 \eqno (5.4)$$
 $V_{d_1,...,d_{n+1}}$ is a hypersurface in ${\bf CP}^{n+1}$ with $D^n$
 isolated singularities each of which is equivalent to the singularity
 at the origin of $q(x_1,...,x_{n+1})=x_1 ^{d_1}+...x_{n+1}^{d_{n+1}}$.
 For a generic hyperplane $H \subset {\bf CP}^{n+1}$ one has
the isomorphism:
 $$\pi_n ({\bf CP}^{n+1}-(V \cup H)) \otimes {\bf Q} =
H_n (M_q, {\bf Q}) \eqno (5.5)$$ where $M_q$ is the Milnor fiber of the
singularity of $q$ at the origin.
 This isomorphism is the isomorphism of ${\bf Q}
[t,t^{-1}]$-modules where the module structure on the right is the one
in which $t$ acts  as the monodromy operator.
\bigskip
\par {\bf Proof}. The hypersurface $V_{{d_1},...,{d_{n+1}}}$ is a
section by a generic linear subspace of dimension $n+1$
of a hypesurface $\tilde V_{d_1,...d_{n+1}}$
in ${\bf CP}^{(n+2)(n+1)-1}$ given by:
 $$G_{d_1,...,d_{n+1}}=
 \tilde g_{d_2 \cdot ... \cdot d_{n+1}}^{d_1}(z_{1,0},z_{1,1},...
 ,z_{1,n+1})+....+\tilde g_{d_1 \cdot ...\cdot d_{n}}^{d_{n+1}}
 (z_{n+1,0},...,z_{n+1,n+1})=0 \eqno (5.6)$$
 where $\tilde g_k$ are generic forms of disjoint set of variables.
 The hypersurface (5.6) has a singular locus containing the component
 of codimension $n+1$ in the ambient space which is given by
  $\tilde g_1=... \tilde g_{n+1}=0$ (i.e. having codimension  $n$ inside
 the hypersurface) as well as possibly some
  components of larger codimensions.
 According to lemma  4.9 for generic hyperplane $\tilde H$
 the module $\pi _i ({\bf CP}^{(n+2)(n+1)-1}
 -(\tilde V_{d_1,...,d_{n+1}} \cup \tilde H)) \otimes {\bf Q} $
 is isomorphic to  $H_n (M_{G_{d_1,...,d_{n+1}}},{\bf Q})$ with the
 usual ${\bf Q} [t,t^{-1}]$-module structure. The forms
 $\tilde g_{d_2 \cdot ...\cdot d_{n+1}},...,
 \tilde g_{d_1 \cdot ... \cdot d_n}$ have isolated singularities at the
 origins of corresponding  ${\bf C}^{n+2}$ because of the genericity
 assumption. Hence, the preceeding lemma implies that $H_n (M_{G_{d_1,...
 ,d_{n+1}}},{\bf Q})$ is isomorphic to $H_n (M_q,{\bf Q})$.
\bigskip
\par {\bf Proposition 5.3.}  Let $f_i=0$ ($i=1,2$) be the equation
 of a hypersurface $V_{f_i}$ of a degree $d$
 in ${\bf CP}^{n_i+1}$. Assume that
 the codimension of the singular locus of $V_{f_i}$ is $k_i$. Then
 the hypersurface $V_{f_1+f_2}$ in ${\bf CP}^{n_1+n_2+3}$ given by
 $f_1+f_2=0$ has the singular locus of codimension $k_1+k_2+1$ and
  $$\pi _{k_1+k_2+1} ({\bf CP}^{n_1+n_2+3}-V_{f_1+f_2})
 \otimes {\bf Q} =(\pi_{k_1} ({\bf CP}^{n_1+1}-V_{f_1})  \otimes {\bf Q})
 \otimes _ {\bf Q} (\pi_{k_2} ({\bf CP}^{n_2+1}-V_{f_2}) \otimes {\bf Q})
   \eqno (5.7)$$.
\bigskip
\par {\bf Proof}. This is an immediate consequence of (4.9) and the
 Thom-Sebastiani theorem.
\bigskip
\par {\bf Examples 5.4.} 1. Let $f(x_0,x_1,x_2)=0$ be an
equation of a curve $C$ of degree 6  which has 6 cusps on a conic. The
 homology of the infinite cyclic cover of ${\bf CP}^2- (C \cup L)$ for
 generic line is ${\bf Q} [t,t^{-1}]/(t^2-t+1)$. (cf. [L]. Recall
 that $\pi _1 ({\bf CP}^2 -(C \cup L))$ is the braid group on
 3 strings, i.e the group of the trefoil knot, and hence the homology
 in question is the Alexander module of the  trefoil. Let $g(y_0,
 y_1,y_2)= 0$ be an equation of another sextic with six cusps
  also on a conic. According to proposition (5.3)
  (in which the homotopy groups in case of curves
 are replaced by the Alexander modules) the generic section by
 ${\bf CP}^4$ of the hypersurface in ${\bf CP}^6$ given by
 $$f(x_0,x_1,x_2)+g(y_0,y_1,y_2,)=0 \eqno (5.8)$$
  is a threefold $W$ which  has isolated singularities (the number of
   which is $6^3$) and  the order of the homotopy group  $\pi_3 \otimes
 {\bf Q}$ of the complement is $(t^2-t+1)$. If one takes as $f$ in (5.8)
 the equation of a sextic with nine cusps which is dual to
  a non singular cubic and uses the fact that
 the Alexander module for its complement is $({\bf Q} [t,t^{-1}]/
(t^2-t+1))^ {\oplus 3}$. This is a consequence of the calculation of the
fundamental group of the complement to such a curve due to O.Zariski.
 He found that the fundamental  group of this curve
 is the kernel of the canonical map of the
braid group of the torus $S^1 \times S^1$ onto $H_1 (S^1 \times S^1,
{\bf Z})$ (cf. [Z]) Combining this with  the calculation which
 uses the  Fox calculus and the
presentation for the braid group of the torus one arrives to
 the Alexander module as above.
 Therefore  one obtains the threefold $W$ with $\pi_3 ({\bf CP}^4-W)
 \otimes {\bf Q}$ being the same as the Alexander module just mentioned.
 Iteration of this
  construction obtained by replacing in  (5.8) $g$ by the equation
 of an $n$-dimensional hypersurface with isolated singularities
  and non vanishing homotopy group $\pi_n \otimes {\bf Q}$
  of the  complement gives examples of hypersurfaces of degree $6$
 and dimension $n$  for which $\pi_n \otimes {\bf Q}$
  has arbitrary large rank for sufficiently large $n$.
 \par 2. Let us consider the equation (5.8) in which $f(x_0,x_1,x_2)$
 is the form giving the equation of a sextic with six cusps not on conic.
 The fundamental group of the complement to such a curve is abelian (cf.
 [Z]) and therefore the homology of the universal cyclic cover of the
  complement is trivial. Hence in this case the construction of example 1
 yields a threefold $W'$ which has $\pi_3 ({\bf CP}^4-W')
 \otimes {\bf Q}=0$ but which has
 the same degree ( i.e. $6$)  and the same number of singularities
 (i.e. $36$) of the same type as $W$
  (i.e. locally given by $x^2+y^2+u^3+v^3=0$).
\bigskip
\bigskip

\ssectitle{\  References.}

\bigskip
\par [Ab] {\bf S.Abhyankar},Tame coverings and
 fundamental groups of algebraic
 varieties, I, II, III. Amer.J. of Math. {\bf 81}(1959),46-94,{\bf 82},
 120-178,179-190.
\par [A] {\bf M.Atiyah},(notes by S.Donaldson)
Representations of the braid
 group. London Math. Soc. Lecture  Notes Series,vol.150 S.Donaldson
 and C. Thomas ed.115-120 (1990).
\par [ABK] {\bf P.Antonelli, D.Burghelea, P.Kahn.}
 The non-finite homotopy
type of some diffeomorphism groups, Topology.vol.11, No.1,1972.
\par [BK] {\bf E.Brieskorn,H.Knorer},Plane Algebraic
 Curves,Birkhauser-Verlag
 ,Basel-Boston-Stuttgart,1986.
\par [De] {\bf P.Deligne}, Le groupe fundamentale du complementaire d'une
 courbe plane n'ayant que des points doubles ordinaire est abelien,
, Seminar Bourbaki. 1979/1980. Lect.Notes in Math. vol.842, Springer
Verlag,1981,pp1-10.
\par [Deg] {\bf A.Degtyarev},Alexander polynomial of an algebraic
 hypesurface, (in russian) ,Preprint, Leningrad, 1986.
\par [Di] {\bf A.Dimca},Alexander polynomials for
 projective hypersurfaces,
Preprint, Max Plank Institute, Bonn,1991.
\par [Dy] {\bf M.Dyer}
,Trees of homotopy types of $(\pi,m)$ complexes. London
 Math.Soc. Lecture Notes Series, 36.p251-255.
\par [F] {\bf W.Fulton}, On the fundamental group of the complement of
 a nodal curve, Ann.of Math. vol. 111, 1980,407-409.
\par [H] {\bf H.Hamm}, Lefschetz theorems
for singular varieties, Proc.Symp.
 in Pure Math. vol.40, part 1, p547-558. (1983).
\par [vK] {\bf E.R.van Kampen}, On the fundamental group
 of an algebraic curve
, Amer.J.of Mathematics,vol.33 (1935).
\par [KW] {\bf R.Kulkarni,J.Wood}, Topology of non-singular complex
 hypersurfaces, Adv.in Math. 35, (1980) 239-263.
 \par [L1] {\bf A.Libgober}
 , Homotopy groups of the complement to singular
 hypersurfaces, Bull. Amer. Math. Soc. vol. 13, No.1, 1986, p.49-51.
\par [L2] {\bf A.Libgober}
, Alexander polynomials of plane algebraic curves
 and cyclic multiple planes,Duke Math.Journ. vol.49 (1982).
\par [L3] {\bf A.Libgober} On $\pi_2$ of the complements to
 hypersurfaces which are generic projections. Adv. Study in Pure
 Mathematics 8, 1986, Complex Analytic Singularities, pp.229-240.
\par [Mi] {\bf J.Milnor},
 Singular points of complex hypersurfaces, Princeton
 University Press, 1968.
\par [Mo] {\bf B.Moishezon}
,Stable branch curves and braid monodromies,Lect.
 Notes in Mathematics,vol.862.
\par [N] {\bf  M.Nori}
,Zariski conjecture and related problems, Ann.Sci.Ecole
Normal Sup.tome XVI, p.305-344.
\par [R]{\bf  R.Randell},Milnor fibres and Alexander polynomials of plane
curves, Proc.Symp.Pure Math. vol.40,p.415-420. (1983).
\par [Wh] {\bf J.H.C.Whitehead}
,Simple homotopy type, Amer.J. of Mathematics,
 72 (1950),1-57.
\par [Z] {\bf O.Zariski}, Chapter 8 in Algebraic Surfaces, Sprinder
Verlag, 1972.
   \end
 \end